\documentclass[final,5p,times,twocolumn]{elsarticle}

\usepackage{amsmath,amssymb,amsfonts}
\usepackage{graphicx}          
\usepackage{color}
\usepackage{thmtools}

\usepackage{lipsum}

\declaretheorem{example}
\renewcommand\thmcontinues[1]{Continued} 

\newtheorem{thm}{Theorem}
\newtheorem{lem}{Lemma}
\newproof{pf}{Proof}
\newtheorem{prop}{Proposition}
\newtheorem{remark}{Remark}

\newcommand{\tm}{\times}
\newcommand{\trn}{^{\scriptscriptstyle \top}}
\newcommand{\N}{\mathbb{N}}
\newcommand{\Z}{\mathbb{Z}}
\newcommand{\R}{\mathbb{R}}

\newcommand{\ep}{\varepsilon}

\newcommand{\AC}{\mathcal{A}}
\newcommand{\BC}{\mathcal{B}}
\newcommand{\CC}{\mathcal{C}}
\newcommand{\GC}{\mathcal{G}}
\newcommand{\JC}{\mathcal{J}}
\newcommand{\NC}{\mathcal{N}}

\newcommand{\SC}{\mathcal{S}}
\newcommand{\TC}{\mathcal{T}}
\newcommand{\WC}{\mathcal{W}}

\newcommand{\inv}{\mathrm{inv}}

\newcommand{\colorb}{}

\usepackage{hyperref}

\usepackage{algorithm}
\usepackage{algpseudocode}
\usepackage{caption}

\begin{document}

\begin{frontmatter}

\title{Numerical over-approximation of invariance entropy via finite abstractions}

\cortext[cor1]{This work is supported in part by the NSF
under Grant CMMI-2013969 and by the German Research Foundation (DFG) through grants RU 2229/1-1 and ZA 873/4-1. Corresponding author: M.~S.~Tomar}

\author[bou]{M.~S.~Tomar}\ead{mahendra.tomar@colorado.edu}    
\author[lmu]{C.~Kawan}\ead{christoph.kawan@lmu.de}               
\author[bou,lmu]{M.~Zamani}\ead{majid.zamani@colorado.edu}  

\address[bou]{Computer Science Department, University of Colorado Boulder, USA}
\address[lmu]{Institute of Informatics, LMU Munich, Germany}
          

\begin{abstract}                          
For a closed-loop system with a digital channel between the sensor and controller, invariance entropy quantifies the smallest average rate of information above which a compact subset $Q$ of the state set can be made invariant. There exist different versions of invariance entropy for deterministic and uncertain control systems, which are equivalent in the deterministic case. In this paper, we present the first numerical approaches to obtain rigorous upper bounds of these quantities. Our approaches are based on set-valued numerical analysis and graph-theoretic constructions. We combine existing algorithms from the literature to carry out our computations for several linear and nonlinear examples. A comparison with the theoretical values of the entropy shows that our bounds are of the same order of magnitude as the actual values.
\end{abstract}

\end{frontmatter}

\section{Introduction}

In classical control theory, sensors and controllers are usually connected through point-to-point wiring. In networked control systems (NCS), sensors and controllers are often spatially distributed and involve wireless digital communication networks for data transfer. Compared to classical control systems, NCS provide many advantages such as reduced wiring, low installation and maintenance costs,
greater system flexibility and ease of modification. NCS find applications in many areas such as car automation, intelligent buildings, and transportation networks. However, the use of communication networks in feedback control loops makes the analysis and design of NCS much more complex. In NCS, the use of digital channels for data transfer from sensors to controllers limits the amount of data that can be
transferred per unit of time. This introduces quantization errors that can affect the control performance adversely.%

The problem of stabilizing or observing a system over a communication channel with a limited bit rate has attracted a lot of attentions in the past two decades. In this context, a classical result, often called the \emph{data-rate theorem}, states that the minimal bit rate or channel capacity above which a linear system can be stabilized or observed is given by the log-sum of the unstable eigenvalues. This result has been proved under various assumptions on the system model, channel model, communication protocol, and stabilization/estimation objectives. Comprehensive reviews of results on data-rate-limited control can be found, e.g., in the surveys \cite{andrievsky2010control,franceschetti2014elements,nair2007feedback} and the books \cite{yuksel2013stochastic,matveev2009estimation}.%

For nonlinear systems, the smallest bit rate of a channel between the coder and the controller, to achieve some control task such as stabilization or invariance, can be characterized in terms of certain notions of \emph{entropy} which are defined in terms of the open-loop system and are independent of the choice of the coder-controller. In spirit, they are similar to classical entropy notions used in the theory of dynamical systems to quantify the rate at which a system generates information, see \cite{katok2007fifty}.%

In this paper, we first consider deterministic systems and focus on the notion of \emph{invariance entropy (IE)} introduced in \cite{ColoniusKawan09} as a measure for the smallest average data rate above which a compact controlled invariant subset $Q$ of the state set can be made invariant. We present the first attempt to compute upper bounds of IE numerically. Our approach combines different algorithms. First, we compute a symbolic abstraction of the given control system over the set $Q$ and the corresponding invariant controller using the tool \texttt{SCOTS}~\cite{rungger2016scots}. This results in a fine box partition of $Q$ with a set of admissible control inputs assigned to each box for maintaining invariance of $Q$. In the second step, we use the tool \texttt{dtControl}~\cite{ashok2020dtcontrol} that converts the controller from a look-up table into a decision tree. Each leaf node of the tree represents a group of boxes to which the same single control input is assigned. The set of groups constitute a coarse partition of $Q$. Finally, in the third step, an algorithm that was proposed in~\cite{froyland2001rigorous} for estimation of topological entropy is adopted. Its output serves as an upper bound for the IE.%

In addition, we also develop a method to approximate the IE of uncertain control systems, as introduced in~\cite{rungger2017invariance,tomar2020invariance}, that generalizes the IE of deterministic systems. If the IE of a set $Q$ (for an uncertain system) is finite~\cite[Sec. 4.2]{rungger2017invariance}, an upper bound can be computed from a graph constructed using a finite abstraction of the system~\cite[Sec.~6]{rungger2017invariance}. However, the number of vertices in the graph is of the order of $2^{2^n}$, where $n$ is the number of states in the finite abstraction. In this paper, we present an upper bound for the IE of uncertain systems that can be computed from a weighted directed graph constructed from an invariant partition (a pair of a finite partition of $Q$ and a map that assigns a control input to every partition element). Our main result characterizes the entropy of the invariant partition in terms of the weights of the graph and establishes that it is the same as the maximum cycle mean of the graph. We should highlight that the number of vertices in this graph is not larger than $n$. Our proposed procedures may still suffer from the curse of dimensionality due to constructing finite abstractions of control systems. Moreover, at this point, we are not able to quantify the gap between the upper bounds and the actual values of the IE.%

{\bf Brief literature review.} The notion of invariance entropy for deterministic systems is equivalent to \emph{topological feedback entropy} that has been introduced earlier in \cite{NairEvansMarrelsMoran04}; see \cite{colonius2013note} for a proof. Various notions of invariance entropy have been proposed to tackle different control problems or other classes of systems, see for instance \cite{Colonius12minimal} (exponential stabilization), \cite{kawan2015network} (invariance in networks of systems), \cite{rungger2017invariance} (invariance for uncertain systems), \cite{colonius2018metric,wang2019measure} (measure-theoretic versions of invariance entropy) and \cite{kawan2019invariance} (stochastic stabilization). An over-approximation of invariance entropy through a compositional approach, for networks of uncertain control systems, was also discussed in \cite{tomar2020compositional}. Also the problem of state estimation over digital channels has been studied extensively by several groups of researchers. As it turns out, the classical notions of entropy used in dynamical systems, namely measure-theoretic and topological entropy (or variations of them), can be used to describe the smallest data rate or channel capacity above which the state of an autonomous dynamical system can be estimated with an arbitrarily small error, see \cite{savkin2006analysis,liberzon2017entropy,sibai2017optimal,yang2018topological,kawan2018optimal}. Motivated by the observation that estimation schemes based on topological entropy suffer from a lack of robustness and are hard to implement, the authors of \cite{MaPo_automatica,PartII} introduce a suitable notion of \emph{restoration entropy} which characterizes the minimal data rate for so-called \emph{regular} and \emph{fine observability}. Finally, algorithms for state estimation over digital channels have been proposed in several works, particularly in \cite{liberzon2017entropy,MaPo_automatica,hafstein2019numerical,kawan2021subgradient}.%

\textbf{Related work:} In~\cite{gao2021invariant}, the authors consider linear uncertain control systems and provide an algorithm to compute an invariant cover, the cardinality of which serves as an upper bound for the invariance entropy. In contrast, our proposed procedure here is applicable to nonlinear systems as well.%

{\bf Notation:} We write $\N = \{1,2,3,\ldots\}$ for the natural numbers, $\Z$ for the set of integers and $\Z_+ := \N \cup \{0\}$. By $\R$, we denote the set of real numbers and define $\R_+ := \{ r\in \R: r\geq 0\}$ and $\R_{>0} := \R_+ \backslash \{0\}$. By $[a;b] = \Z \cap [a,b]$ and $[a;b) = \Z \cap [a,b)$, we denote closed and right-open discrete intervals. We write $|A|$ for the cardinality of a set $A$ and $\rho(R)$ for the spectral radius of a square matrix $R$. The notation $Y^X$ is used for the set of all functions $f:X\rightarrow Y$. {For $\tau\in\Z_+$, we use $X^\tau$ to denote $X^{[0;\tau)}$.} By $f:X\rightrightarrows Y$, we denote a set-valued map from $X$ to $Y$.
{A \emph{cover} $\AC$ of set $Q$ is a family of subsets of $Q$ such that $\cup_{A\in \AC} A = Q$. A cover $\AC$ is called a \emph{partition} if for all $A_1, A_2\in {\cal A}$, $A_i\neq \emptyset$ and $A_1\cap A_2 = \emptyset$. We write $f|_M$ for the restriction of a map $f$ to a subset $M\subseteq X$.}

\section{Background on invariance entropy}\label{sec_background}

In this section, we provide the necessary theoretical background for our proposed numerical methods.%

A \emph{deterministic discrete-time control system} is given by%
\begin{equation}\label{eq_det_sys}
  \Sigma:\quad x_{t+1} = f(x_t,u_t),%
\end{equation}
where $f:X \tm U \rightarrow X$, $X \subseteq \R^n$, $U \subseteq \R^m$, is a (not necessarily continuous) map. The \emph{transition map} $\varphi:\Z_+ \tm X \tm U^{\Z_+} \rightarrow X$ of $\Sigma$ is defined as%
\begin{equation*}
  \varphi(t,x,\omega) := \left\{\begin{array}{cl}
	                             x & \mbox{if } t = 0,\\
															f(\varphi(t-1,x,\omega),\omega_{t-1}) & \mbox{if } t > 0.
														\end{array}\right.%
\end{equation*}
Now, consider a compact set $Q \subseteq X$ which is \emph{controlled invariant}, i.e., for each $x \in Q$ there is $u \in U$ with $f(x,u) \in Q$. For any $\tau \in \N$, a set $\SC \subset U^{\tau}$ is called \emph{$(\tau,Q)$-spanning} if for each $x \in Q$ there is $\omega \in \SC$ with $\varphi(t,x,\omega) \in Q$ for $0 \leq t \leq \tau$. We write $r_{\inv}(\tau,Q)$ for the minimal cardinality among all $(\tau,Q)$-spanning sets and define the \emph{invariance entropy (IE)} of $Q$ as%
\begin{equation*}
  h_{\inv}(Q) := \lim_{\tau \rightarrow \infty}\frac{1}{\tau}\log_2 r_{\inv}(\tau,Q),%
\end{equation*}
if $r_{\inv}(\tau,Q)$ is finite for all $\tau$; otherwise, $h_{\inv}(Q) := \infty$. The existence of the limit follows from the subadditivity of the sequence $(\log_2 r_{\inv}(\tau,Q))_{\tau\in\N}$, using Fekete's subadditivity lemma (see \cite[Lem.~2.1]{colonius2013note} for a proof).%

The method we propose to estimate $h_{\inv}(Q)$ is based on an alternative characterization of this quantity that we will now describe. A triple $(\AC,\tau,G)$ is called an \emph{invariant partition} of $Q$ if $\AC$ is a finite partition of $Q$, $\tau \in \N$, and $G:\AC \rightarrow U^{\tau}$ is a map satisfying\footnote{\colorb For a set $A\subset X$, by $\varphi(t,A,G(A))$ we refer to $\cup_{x\in A} \varphi(t,x,G(A))$. } $\varphi(t,A,G(A))\subseteq Q$ for every $A \in \AC$ and $0 \leq t \leq \tau$ (note that $\varphi(t,x,\omega)$ only depends on $\omega|_{[0;t)}$). For a given $\CC = (\AC,\tau,G)$, we define%
\begin{equation*}
  T_{\CC}:Q \to Q,\quad T_{\CC}(x) := \varphi(\tau,x,G(A_x)),%
\end{equation*}
where $A_x \in \AC$ is such that $x \in A_x$. Since $\AC$ is a partition of $Q$, $T_{\CC}$ is well-defined.%

Now, let $\CC = (\AC,\tau,G)$ be an invariant partition. For each $N \in \N$, we introduce the set%
\begin{align*}
  {\colorb \WC_N(T_\CC)} &:= \{ \alpha \in \AC^N : \exists x \in Q \\
	 &\qquad\qquad\qquad \mbox{s.t.\ } T_{\CC}^i(x) \in \alpha_i,\ 0 \leq i < N\},
\end{align*}
{\colorb which is constituted by all such $N$-length sequences in $\AC$ that there exists a trajectory of $T_{\CC}$ that follows the sequence.} Next we define%
\begin{equation*}
  {\colorb h^*(T_\CC) } := \lim_{N \rightarrow \infty}\frac{1}{N}\log_2 |\WC_N(T_\CC)|.%
\end{equation*}
Again, subadditivity guarantees the existence of the limit. Then, by \cite[Thm.~2.3]{kawan2013invariance}, the IE satisfies%
\begin{equation}\label{eq_hinv_alt}
  h_{\inv}(Q) = \inf_{\CC = (\AC,\tau,G)} \frac{1}{\tau} h^*(T_{\CC}),%
\end{equation}
where the infimum is taken over all invariant partitions of $Q$. In particular, $h_{\inv}(Q) < \infty$ if and only if an invariant partition of $Q$ exists \cite[Prop.~2.20, Lem.~2.3]{kawan2013invariance}.%

An \emph{uncertain discrete-time control system} is given by%
\begin{equation}\label{eq_unc_sys}
  \Sigma:\quad x_{t+1} \in F(x_t,u_t),%
\end{equation}
where $X \subseteq \R^n$, $U \subseteq \R^m$, and $F:X \tm U \rightrightarrows X$ is a set-valued map satisfying $F(x,u) \neq \emptyset$ for all $(x,u) \in X \tm U$.%

Consider a compact set $Q \subseteq X$ which is controlled invariant, i.e., for each $x \in Q$ there is $u \in U$ with $F(x,u) \subseteq Q$. We define the invariance entropy of $Q$ in a quite different manner as in the deterministic case. However, in the special case when $F$ is single-valued, i.e., when $\Sigma$ is deterministic, it coincides with the previous notion.%

A pair $(\AC,G)$ is called an \emph{invariant cover} of $Q$ (w.r.t.~$\Sigma$) if $\AC$ is a finite cover of $Q$ and $F(A,G(A)) \subseteq Q$ for all $A \in \AC$. In the case when $\AC$ is a partition, we call $(\AC,G)$ an \emph{invariant partition}, analogously to the deterministic case.\footnote{However, for uncertain systems, time steps larger than $1$ should be avoided, so there is no number $\tau$ here.} For $\tau \in \N$, let $\JC \subseteq \AC^{[0;\tau)}$ be a set of sequences in $\AC$ of length $\tau$. For $\alpha \in \JC$ and $t \in [0;\tau-2]$, define%
\begin{align}\label{eq_def_pj_1}
\begin{split}
  P_{\JC}(\alpha|_{[0;t]}) &:= \{ A \in \AC : \alpha|_{[0;t]}A = \hat{\alpha}|_{[0;t+1]} \\
	&\qquad\qquad\qquad\qquad \mbox{for some } \hat{\alpha} \in \JC \},%
\end{split}
\end{align}
as the set of immediate successor cover elements $A$ of $\alpha|_{[0;t]}$ in $\JC$, and for $t = \tau-1$, define%
\begin{align*}
  P_{\JC}(\alpha|_{[0;t]}) = P_{\JC}(\alpha) := \{ A \in \AC : A &= \hat{\alpha}(0) \\
	& \mbox{for some } \hat{\alpha} \in \JC \},%
\end{align*}
as the set of the first components of the sequences in $\JC$. Although this set does not depend on $\alpha$, for consistency reasons, we still use the same notation as in \eqref{eq_def_pj_1}. A set $\JC \subseteq \AC^{[0;\tau)}$ is called \emph{$(\tau,Q)$-spanning in $(\AC,G)$} if $P_{\JC}(\alpha)$ covers $Q$ and for all $\alpha \in \JC$ and all $t \in [0;\tau-2]$%
\begin{equation}\label{eq_spanningset_cond}
  F(\alpha(t),G(\alpha(t))) \subseteq \bigcup_{A' \in P_{\JC}(\alpha|_{[0;t]})}A'.%
\end{equation}
In this case, we associate to $\JC$ its \emph{expansion number}%
\begin{equation}\label{eq_def_expno}
  \NC(\JC) := \max_{\alpha\in\JC}\prod_{t=0}^{\tau-1}|P_{\JC}(\alpha|_{[0;t]})|,%
\end{equation}
and write $\bar{r}_{\inv}(\tau,Q,\AC,G)$ for the smallest expansion number among all $(\tau,Q)$-spanning sets in $(\AC,G)$, {\colorb i.e., $\bar{r}_{\inv}(\tau,Q,\AC,G) := \min\{ \NC(\JC)\mid \JC \mbox{ is }(\tau,Q)\mbox{-spanning in }(\AC,G)\} $}. The \emph{entropy of an invariant cover} $(\AC,G)$ is then defined as%
\begin{equation*}
  \bar{h}(\AC,G) := \lim_{\tau \rightarrow \infty}\frac{1}{\tau}\log_2 \bar{r}_{\inv}(\tau,Q,\AC,G).%
\end{equation*}
The existence of the limit follows again by subadditivity. The invariance entropy of $Q$ is now defined as%
\begin{equation*}
  \bar{h}_{\inv}(Q) := \inf_{(\AC,G)}\bar{h}(\AC,G),%
\end{equation*}
where the infimum is taken over all invariant covers of $Q$. Although this definition does not seem to have much similarity with the definition(s) for deterministic systems, $\bar{h}_{\inv}(Q)$ reduces to $h_{\inv}(Q)$ in the case when $F$ is single-valued, see \cite[Thm.~4]{tomar2020invariance}.%

\section{Upper bounds: deterministic case}\label{sec_det_ub}

In this section, we explain how to obtain a computable upper bound for $h_{\inv}(Q)$, based on \eqref{eq_hinv_alt}. Suppose that we have an invariant partition $\CC = (\AC,\tau,G)$ with $\AC = \{A_1,\ldots,A_q\}$ at our disposal. Then any upper bound on $h^*(T_{\CC})$ will yield an upper bound on $h_{\inv}(Q)$.%

Let us first select a refinement $\BC = \{B_1,\ldots,B_r\}$ of $\AC$, i.e., a partition of $Q$ such that each $B\in\BC$ is contained in some $A\in\AC$.
Now we define%
{\colorb
\begin{align*}
  \WC_N(\BC,\AC) :=& \{ \alpha \in \AC^N  : \exists \beta \in \BC^N 
	\mbox{with }T_{\CC}(\beta_j) \cap \beta_{j+1} \neq \emptyset \\ &\forall j \in [0;N-2]  \ \mbox{ s.t. } \beta_i \subseteq \alpha_i\ \forall i \in [0;N-1] \}.%
\end{align*}    }
From \cite[Sec.~2.2]{froyland2001rigorous}, we have%
\begin{equation*}
  h(\BC,\AC) := \lim_{N \rightarrow \infty}\frac{\log_2|\WC_N(\BC,\AC)|}{N} \geq h^*(T_{\CC}).%
\end{equation*}
Moreover, assuming compactness of the partition sets and continuity of the map, it can be shown that $h(\BC,\AC)$ converges to $h^*(T_{\CC})$ as the maximal diameter of the elements of $\BC$ tends to zero, see \cite[Thm.~4]{froyland2001rigorous}.\footnote{Of course, such an assumption can, in general, not be satisfied. We expect that the result still holds true if only a	negligibly small amount of the exponential orbit complexity of the closed-loop dynamics is concentrated on the boundaries of the sets $A_i$. If $T_{\CC}$ was continuous on $Q$, this could be formalized by requiring that these boundaries have measure zero w.r.t.~any $T_{\CC}$-invariant Borel probability measure.}


The paper \cite{froyland2001rigorous} describes an algorithm for the exact computation of $h(\BC,\AC)$, based on symbolic dynamics. 
First, we associate a transition matrix to $\BC$ via\footnote{In the language of symbolic dynamics, the matrix $\Gamma$ defines a \emph{subshift of finite type} over the alphabet $\{1,\ldots,r\}$.}%
\begin{equation}\label{eq_def_Gamma}
  \Gamma_{i,j} := \left\{ \begin{array}{cl} 1 & \mbox{if } T_{\CC}(B_i) \cap B_j \neq \emptyset\\
	0 & \mbox{otherwise} \end{array}\right.,\ i,j=1,\ldots,r.%
\end{equation}
Then one constructs a directed labeled graph $\GC$ from the transition matrix $\Gamma$.
The set of nodes is $\BC$ and $\Gamma_{i,j}=1$ indicates that there is a directed edge from $B_i$ to $B_j$. To this edge, we assign the edge label%
\begin{equation}\label{eq_def_edge_labels}
  L(B_i) := j, \mbox{\ where $j$ is such that\ } B_i \subset A_j.%
\end{equation}
Elements of $\WC_N(\BC,\AC)$ are thus generated by concatenating labels along walks of length $N$ on $\GC$.
{\colorb To compute $h(\BC,\AC)$, a \emph{right-resolving graph}\footnote{A labeled graph is right-resolving if, for each vertex, all the outgoing edges have different labels.} $\bar{\GC}$ needs to be determined (see \cite[\S 3.3]{lind2021introduction}), such that the subset of $\mathbb{N}^{\mathbb{Z}}$ generated by concatenation of edge labels along walks in the graph is same for both $\GC$ and $\bar{\GC}$.\footnote{The subset of $\mathbb{N}^{\mathbb{Z}}$ generated by concatenating edge labels along all  walks on $\GC$ forms a \emph{sofic shift} whose topological entropy equals $h(\BC,\AC)$.
}
Each node in the right-resolving graph $\bar{\GC}$  is some subset of $\BC$, while the set of edge labels is identical for both ${\GC}$ and $\bar{\GC}$; for details on its computation, see \cite{froyland2001rigorous}. For any edge in $\bar{\GC}$, with edge label $e$, directed from node $n_1\subset \BC$ to $n_2\subset\BC$, we have $n_2 = \{  B_j \in \BC \mid  B_i \in n_1, L(B_i) = e, \Gamma_{i,j} = 1 \}$, i.e., $n_2$ is the set of those elements of $\BC$ which have an incoming edge from such members of $n_1$ whose image under the map $L$ equals $e$.}
Let $\tilde{n}$ be the number of nodes in $\bar{\GC}$. An associated $\tilde{n}\tm\tilde{n}$ adjacency matrix $R$ is defined as%
\begin{equation*}
  R_{i,j} := \mbox{\#~of edges from node $i$ to node $j$ in $\bar{\GC}$}.%
\end{equation*}
If $\GC$ is strongly connected (i.e., for every pair of nodes $u$ and $v$, there exists a directed path from $u$ to $v$), then by \cite[Prop.~7]{froyland2001rigorous} we have $h(\BC,\AC) = \log_2\rho(R)$. In general, we need to determine the strongly connected components $\GC_1,\ldots,\GC_p$ of $\GC$ and compute a right-resolving graph for each component separately, resulting in adjacency matrices $R^1,\ldots,R^p$. Then (see \cite[Rem.~9]{froyland2001rigorous})%
\begin{equation*}
  h_{\inv}(Q) \leq \frac{1}{\tau} h(\BC,\AC) = \frac{1}{\tau} \max_{1 \leq k \leq p}\log_2\rho(R^k)
\end{equation*}
{\colorb where $p$ is the number of strongly connected components of $\GC$ and $R^k$ is the adjacency matrix of the $k$-th strongly connected component $\GC_k$. }
This leaves us with the problem of constructing an invariant partition $\CC$ with low entropy, in order to obtain an upper bound for $h_{\inv}(Q)$ which is not too conservative. We are not really able to minimize the entropy, but we can produce invariant partitions $(\AC,\tau,G)$ for a given input sequence length $\tau$ with a (potentially) small number of partition elements. Due to the trivial inequality $h^*(T_{\CC}) \leq \log_2|\AC|$, this is desirable.%

We now explain step by step how we determine an invariant partition and compute its entropy. The following example is used to illustrate each step.%

\begin{example}\label{ex1}
Consider the linear control system%
\begin{equation*}
  \Sigma: \quad x_{t+1} = Ax_t + \left[\begin{array}{c} 1 \\ 1 \end{array}\right]u_t,\quad A = \left[\begin{array}{cc} 2 & 0 \\ 0 & \frac{1}{2} \end{array}\right],%
\end{equation*}
with $x_t \in \R^2$ and $u_t \in [-1,1]$. For the compact controlled invariant set $Q = [-1,1] \tm [-2,2]$, see \cite[Ex.~21]{colonius2021controllability}, we intend to compute an upper bound of $h_{\inv}(Q)$.
\end{example}

Given a discrete-time system $\Sigma$ as in \eqref{eq_det_sys} and a compact controlled invariant set $Q \subseteq X$, we proceed according to the following steps.%
		   
\begin{figure}[tp]
	\centering
	\includegraphics[scale=.6]{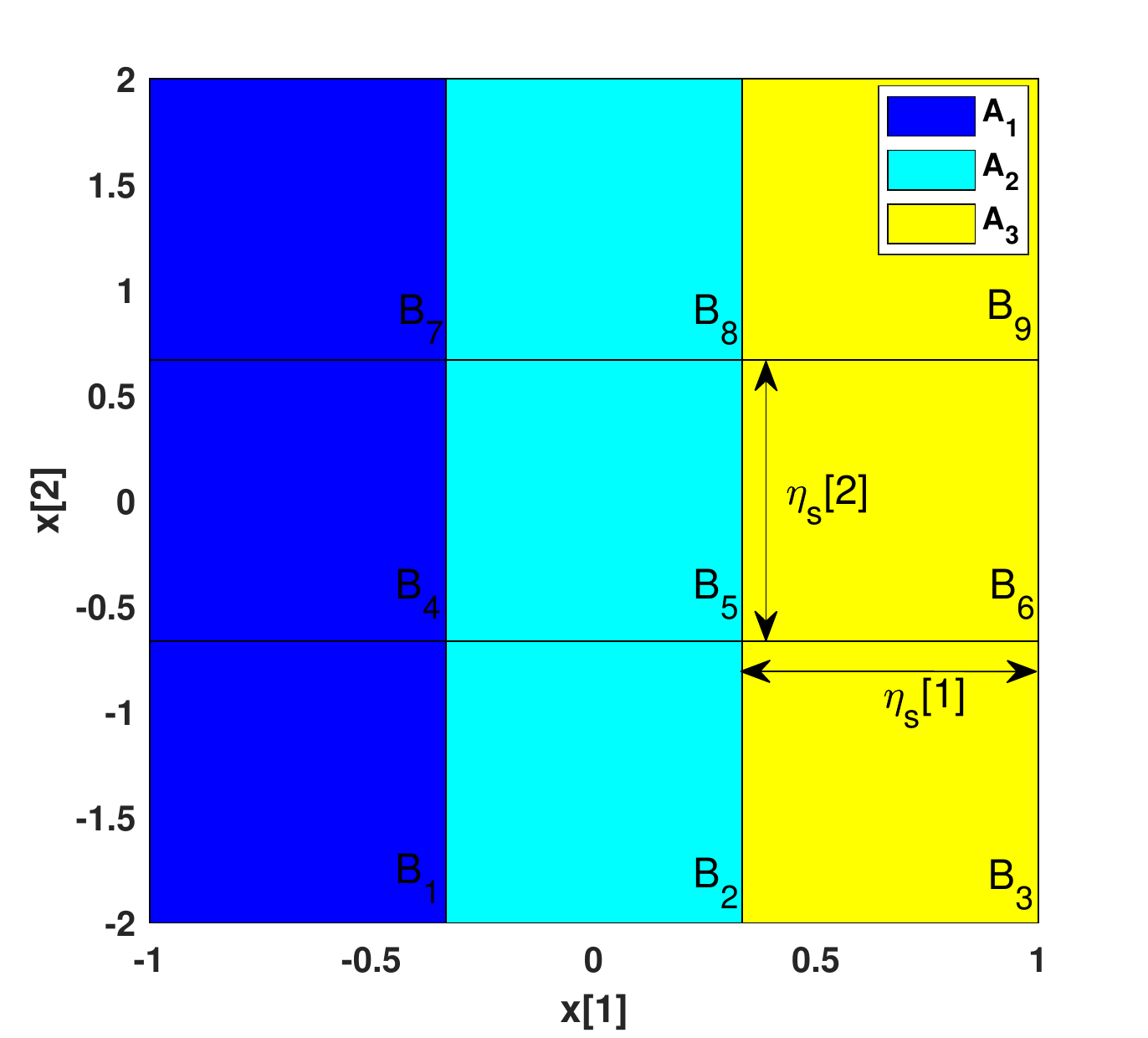}
	\caption{The partitions $\AC$ and $\BC$ for Example \ref{ex1}}
	\label{fig1}
\end{figure}

\begin{enumerate}
\item[(1)] Compute a symbolic invariant controller for the set $Q$. Consider the hyperrectangle $Q_X$ that encloses $Q$ and assume that $Q_X \subseteq X$. We use \verb+SCOTS+ to compute an invariant controller for $\Sigma$ with $Q_X$ as the state set and $\eta_s,\eta_i$ as the grid parameters for the state and input sets, respectively. A smaller value for $\eta_s$ results in a finer grid on the state set, which typically results in a better upper bound. We denote the set of boxes in the domain of the computed controller by $\BC = \{B_1,\ldots,B_r\}$ and put%
\begin{equation*}
  \bar{Q} := \bigcup_{i=1}^r B_i \subseteq Q.%
\end{equation*}
The set $\bar{Q}$ is our approximation of $Q$.%

\begin{example}[continues=ex1]
We used \verb+SCOTS+ with the state set $Q_X = Q$ and the state and input set grid parameters $\eta_s = [2/3,4/3]\trn$ and $\eta_i = 1$. This results in a state set grid with 9 boxes, $\BC = \{B_1,\ldots,B_9\}$ and $\bar{Q} = Q$ (see Fig.~\ref{fig1}).
\end{example}

\item[(2)] The controller obtained in the previous step is, in general, non-deterministic, i.e., different control inputs are assigned to the same state. In this step, we determinize the obtained controller. We denote the closed-loop system ($\Sigma$ with the determinized controller $C$) by $\Sigma_C$. To determinize the controller, we used the state-of-the-art toolbox \verb+dtControl+ \cite{ashok2020dtcontrol}, which utilizes the \emph{decision tree learning algorithm}. This also provides the required coarse partition $\AC$, of which $\BC$ is a refinement.%

\begin{example}[continues=ex1]
For the example, we used \verb+dtControl+ with parameters Classifier $=$ `cart' and Determinizer $=$ `maxfreq'. This results in an invariant partition $(\AC,1,G)$ for the set $\bar{Q} := \bigcup_{B\in\BC}B$, where $\AC$ is a partition of $\bar{Q}$ such that every $A \in \AC$ is a union of some sets in $\BC$ and $G(A) \in U$ is the control input assigned to the set $A$ given by \verb+dtControl+. Figure \ref{fig1} shows the obtained partitions $\AC$ and $\BC$.
\end{example}

\item[(3)] For the dynamical system $\Sigma_C$, we obtain the transition matrix $\Gamma$ (defined in \eqref{eq_def_Gamma}) for the boxes in $\bar{Q}$.%

\item[(4)] We obtain the edge labels map $L(B_i)$ as in \eqref{eq_def_edge_labels}.%

\begin{example}[continues=ex1]
For any $B_i \in \BC$,%
\begin{equation*}
	L(B_i) = \left\{ \begin{array}{cc}
		1 & \mbox{if } i = 1+3t, 0\leq t\leq 2,\\
		2 & \mbox{if } i = 2+3t, 0\leq t\leq 2,\\
		3 & \mbox{if } i = 3+3t, 0\leq t\leq 2.
	\end{array}\right.
\end{equation*}
\end{example}

\item[(5)] We construct a directed labeled graph $\GC$ with $\BC$ as the set of nodes. If $\Gamma_{i,j} = 1$, there is a directed edge from the node $B_i$ to $B_j$ with label $L(B_i)$.%

\item[(6)] We determine the strongly connected components of $\GC$.%

\begin{figure*}
\begin{minipage}{0.5\textwidth}
\centering
	\includegraphics[scale=.3]{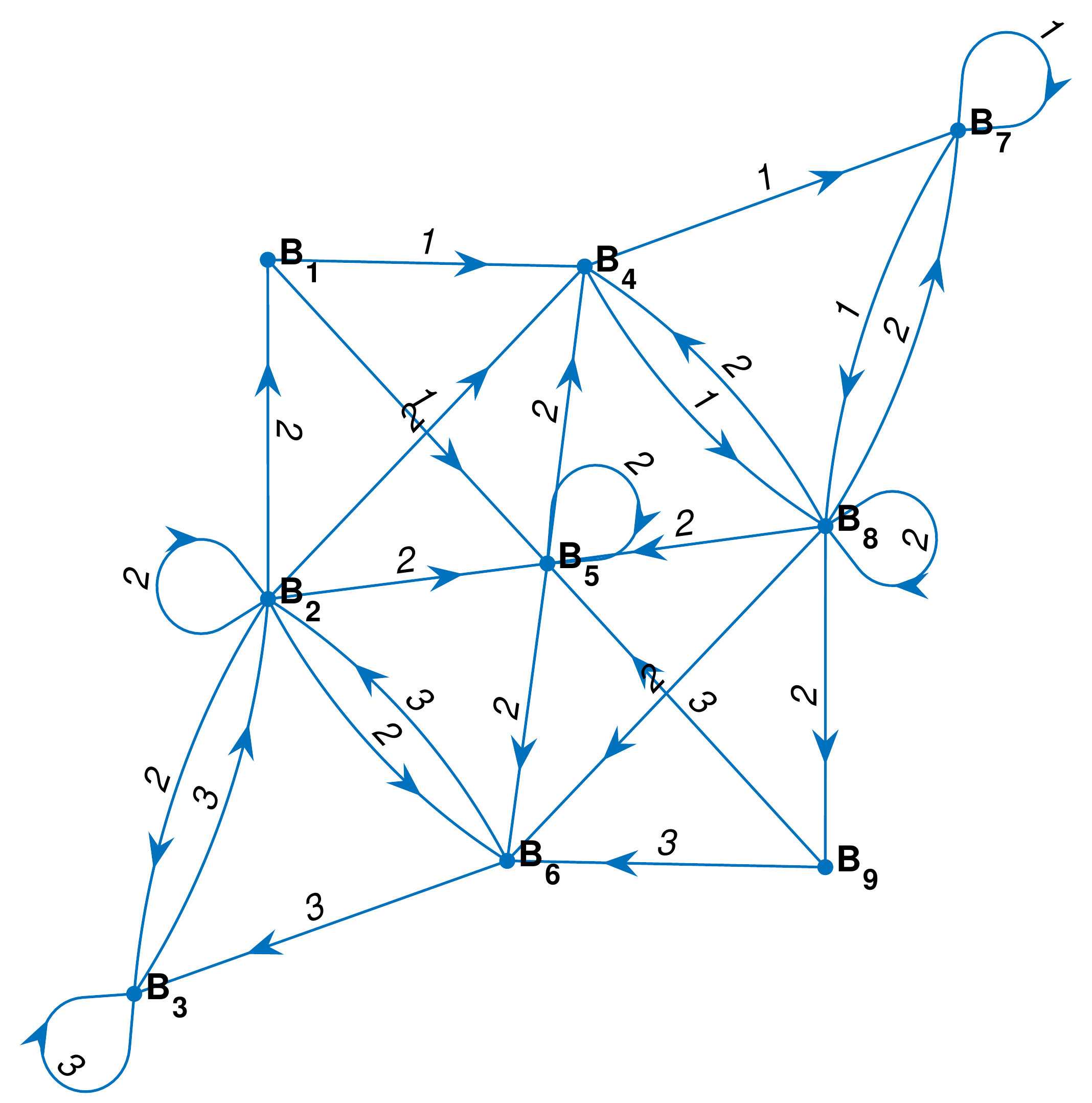}
	\vspace{-5mm}
	\caption{The graph $\GC$ for Example \ref{ex1}}
	\label{fig:G}
\end{minipage}
\begin{minipage}{0.5\textwidth}
	\centering
	\includegraphics[scale=.3]{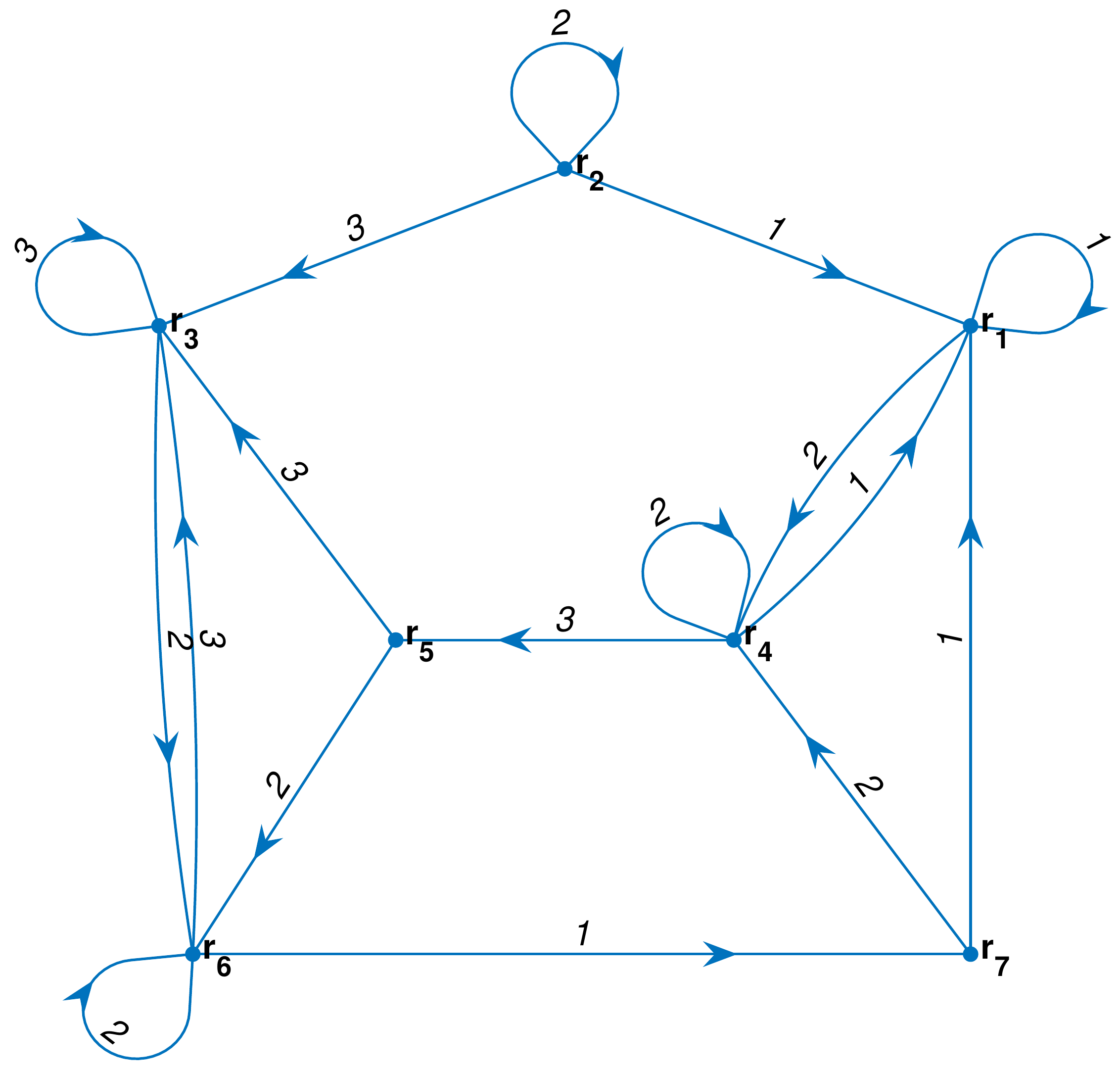}
	\vspace{-5mm}
	\caption{The right-resolving graph $\bar{\GC}$ for Example \ref{ex1}}
	\label{fig:Gbar}
\end{minipage}
\end{figure*}

\begin{example}[continues=ex1]
	$\GC$ is strongly connected. {\colorb Figure~\ref{fig:G} shows the constructed graph $\GC$.}
\end{example}

\item[(7)] For every strongly connected component $\GC_k$ of $\GC$, we find a right-resolving graph $\bar{\GC}_k$. The directed graph $\bar{\GC}_k$ is deterministic in the sense that for every node no two outgoing edges have the same label.%

\begin{example}[continues=ex1]
Right-resolving graph of $\GC$ with nodes $r_1 = \{B_i : i\in\{7,8\}\}$, $r_2 = \{B_i : 4\leq i\leq 6\}$, $r_3 = \{B_i : i\in\{2,3\}\}$, $r_4 = \{B_i : 4\leq i\leq9\}$, $r_5 = \{B_i : i\in\{2,3,5,6\}\}$, $r_6 = \{B_i : 1\leq i\leq 6\}$, and $r_7 = \{B_i : i\in\{4,5,7,8\}\}$. {\colorb The constructed right-resolving graph $\bar{\GC}$ is shown in Figure~\ref{fig:Gbar}. }
\end{example}

\item[(8)] Using $\bar{\GC}_k$, we construct an adjacency matrix $R^k$ by $R^k_{i,j} := l$, where $l$ is the number of edges from node $i$ to node $j$ in $\bar{\GC}_k$.%
  
\begin{example}[continues=ex1]
From $\bar{\GC}_R$, we obtain%
\begin{equation*}
	R = \left[\begin{array}{ccccccc}
		1 & 0 & 0 & 1 & 0 & 0 & 0 \\
		1 & 1 & 1 & 0 & 0 & 0 & 0 \\
		0 & 0 & 1 & 0 & 0 & 1 & 0 \\
		1 & 0 & 0 & 1 & 1 & 0 & 0 \\
		0 & 0 & 1 & 0 & 0 & 1 & 0 \\
		0 & 0 & 1 & 0 & 0 & 1 & 1 \\
		1 & 0 & 0 & 1 & 0 & 0 & 0 \\
	\end{array}\right],%
\end{equation*}
with $\rho(R) = 2.41421$ and $\log_2(2.4142) = 1.2716$.
\end{example}

\end{enumerate}

\section{Upper bounds: uncertain case}\label{sec_ub_uncertain}

In this section, we explain how to obtain a computable upper bound for the IE of an uncertain system.%

Suppose again that we know an invariant partition $(\AC,G)$ and recall that the time step $\tau$ is always set to $1$ for uncertain systems. We define a set-valued map $T:Q \rightrightarrows Q$ by $T(x) := F(x,G(A_x))$, where $x \in A_x \in \AC$. We also define a weighted directed graph $\GC$ with $\AC$ as its set of nodes. The graph $\GC$ contains an edge from $A$ to $A'$, denoted by $e_{AA'}$, if $T(A) \cap A' \neq \emptyset$. We define maps $D:\AC \rightrightarrows \AC$ and $w:\AC \rightarrow \R_+$ by%
\begin{align}\label{eq_def_weightfnc}
\begin{split}
  D(A) &:= \{A' \in \AC : T(A) \cap A' \neq \emptyset\}, \\
	w(A) &:= \log_2 |D(A)|.%
\end{split}
\end{align}
The weight of the edge $e_{AA'}$ is defined to be $w(A)$. We observe that%
\begin{equation}\label{eq_T_incl}
  T(A) \subseteq \bigcup_{\hat{A} \in D(A)}\hat{A}.%
\end{equation}
Given $\tau \in \N \cup \{\infty\}$, we let $W_{\tau}(\GC)$ denote the set of all (node) paths $(A_i)_{i=0}^{\tau-1}$ in $\GC$ of length $\tau$.%

Consider a cycle $c = (e_{A_iA_{i+1}})_{i=1}^k$, $A_{k+1} = A_1$, in $\GC$. The \emph{mean cycle weight} of $c$ is defined as%
\begin{equation*}
  w_{\rm m}(c) := \frac{1}{k}\sum_{i=1}^k w(A_i).%
\end{equation*}
The \emph{maximum cycle mean weight} is then defined as%
\begin{equation*}
  w^*_{\rm m}(\GC) := \max_c w_{\rm m}(c),%
\end{equation*}
the maximum taken over all cycles in $\GC$ (the maximum exists because, due to the finiteness of the graph, it suffices to take the maximum over finitely many cycles).%

Our algorithm is based on the following theorem, which yields a characterization of the entropy of an invariant partition in terms of the associated graph $\GC$.%

\begin{thm}\label{thm1}
For an uncertain control system $\Sigma$ as in \eqref{eq_unc_sys}, a compact controlled invariant set $Q \subseteq X$ and an invariant partition $(\AC,G)$, we have%
\begin{equation}\label{eq_uie_bounds}
  \bar{h}(\AC,G) = \lim_{\tau \rightarrow \infty}\frac{1}{\tau}\max_{\alpha \in W_{\infty}(\GC)}\sum_{t=0}^{\tau-2}w(\alpha(t)) = w^*_{\rm m}(\GC).%
\end{equation}
\end{thm}

{\colorb
\begin{remark}
In~\cite{ahmadi2012joint}, the authors show that the logarithm of the joint spectral radius of a finite set of rank one matrices equals the maximum cycle mean in a directed complete graph.
The result of the paper can be used to establish that, for the case of non-complete graph, the entropy of an invariant partition is upper bounded by the maximum cycle mean.
\end{remark}    }

The rest of this section is devoted to the proof of the theorem. We start with two lemmas.%

\begin{lem}
$W_{\tau}(\GC)$ is a $(\tau,Q)$-spanning set in $(\AC,G)$.
\end{lem}

\begin{pf}
Since $(\AC,G)$ is an invariant cover, we have $D(A) \neq \emptyset$ for every $A \in \AC$. Thus, for every node in $\GC$ there is at least one outgoing edge. Hence, for all $A\in\AC$ and $\tau\in\N$, there is at least one path of length $\tau$ starting from $A$. It follows that%
\begin{equation*}
  \{ \alpha(0) : \alpha \in W_{\tau}(\GC) \} = \AC.%
\end{equation*}
Consider any $\alpha \in W_{\tau}(\GC)$ and $t \in [0;\tau-1]$. By the definition of $\GC$, we have an edge from $\alpha(t)$ to every $A \in D(\alpha(t))$. Thus, for every $t \in[0;\tau-2]$ we have%
\begin{equation}\label{eq_17}
  P_{W_{\tau}(\GC)}(\alpha|_{[0;t]}) = D(\alpha(t)).%
\end{equation}
Using \eqref{eq_T_incl}, we conclude that $W_{\tau}(\GC)$ satisfies \eqref{eq_spanningset_cond}, and hence is a $(\tau,Q)$-spanning set in $(\AC,G)$. \hfill \qed
\end{pf}

\begin{lem}\label{lem_spanningsets}
For any $(\tau,Q)$-spanning set $\JC$ in $(\AC,G)$, $W_{\tau}(\GC) \subset \JC$.
\end{lem}

\begin{pf}
Let $\JC$ be a $(\tau,Q)$-spanning set in $(\AC,G)$. Then, since $\AC$ is a partition, $\{\alpha(0) : \alpha \in \JC\} = \AC$. If $\alpha \in \JC$ and $t \in [0;\tau-1]$, then from \eqref{eq_spanningset_cond} it follows that $P_{\JC}(\alpha|_{[0;t]})$ covers $F(\alpha(t),G(\alpha(t))) = T(\alpha(t))$. Since $\AC$ is a partition, $D(\alpha(t))$ must be contained in every subset of $\AC$ that covers $T(\alpha(t))$, thus $P_{\JC}(\alpha|_{[0;t]}) \supseteq D(\alpha(t))$. Let $\beta \in W_{\tau}(\GC)$. Then $\beta(0) \in \AC = \{\alpha(0) : \alpha \in \JC\}$, implying $\beta(0) = \alpha(0)$ for some $\alpha \in \JC$. From \eqref{eq_17}, we have $P_{W_{\tau}(\GC)}(\beta(0)) = D(\beta(0))$. Similarly to the reasoning above, since $\AC$ is a partition, $D(\beta(0))$ is contained in every subset of $\AC$ which covers $T(\beta(0))$. As $\JC$ is $(\tau,Q)$-spanning, from \eqref{eq_spanningset_cond} we know that $T(\alpha(0))$ is covered by $P_{\JC}(\alpha(0))$, implying $P_{\JC}(\alpha(0)) \supseteq D(\beta(0))$. From the definition of $\GC$, we obtain $\beta(1) \in D(\beta(0))$, which leads to $\beta(1) \in P_{\JC}(\alpha(0))$. Thus, there exists an $\alpha\in\JC$ with $\alpha|_{[0;1]} = \beta|_{[0;1]}$. Inductively, we obtain the existence of $\alpha \in \JC$ with $\alpha = \beta$, which concludes the proof. \hfill \qed
\end{pf}

We can now prove Theorem \ref{thm1}.%

\begin{pf} (of Theorem \ref{thm1}) 
From \eqref{eq_def_expno} and Lemma \ref{lem_spanningsets}, we conclude that for every $(\tau,Q)$-spanning set $\JC$ in $(\AC,G)$, the inequality $\NC(W_{\tau}(\GC)) \leq \NC(\JC)$ holds, implying that
\begin{equation}\label{eq_thm1_ingr1}
  \bar{r}_{\inv}(\tau,Q,\AC,G) = \NC(W_{\tau}(\GC)) \mbox{\quad for all\ } \tau \in \N.%
\end{equation}
By taking logarithms on both sides of \eqref{eq_def_expno} and using \eqref{eq_17} and \eqref{eq_def_weightfnc}, we obtain%
\begin{equation}\label{eq_thm1_ingr2}
  \log_2 \NC(W_{\tau}(\GC) ) = \max_{\alpha \in W_{\tau}(\GC)}\sum_{t=0}^{\tau-2}w(\alpha(t)) + \log_2|\AC|.%
\end{equation}
Putting \eqref{eq_thm1_ingr1} and \eqref{eq_thm1_ingr2} together, it follows that%
\begin{equation*}
  \bar{h}(\AC,G) = \lim_{\tau \rightarrow \infty} \frac{1}{\tau} \max_{\alpha \in W_{\tau}(\GC)}\sum_{t=0}^{\tau-2}w(\alpha(t)).%
\end{equation*}
Observing that the elements of $W_{\tau}(\GC)$ are restrictions of elements of $W_{\infty}(\GC)$ to $[0;\tau-1]$, the first equality in \eqref{eq_uie_bounds} follows.%

For the proof of the second equality in \eqref{eq_uie_bounds}, let $\AC = \{A_1,\ldots,A_q\}$ and consider an arbitrary $\alpha \in W_{\infty}(\GC)$. From \cite[Lem.~3]{traiger1968asymptotic}, we know that for each $\tau$ we can write%
\begin{equation*}
  \sum_{t=0}^{\tau-2} w(\alpha(t)) = \sum_{t=0}^a w(\beta(t)) + \sum_{i=1}^r l_i w_{\rm m}(\sigma_i),%
\end{equation*}
for some $\beta \in W_{\infty}(\GC)$, $a < n - 1$ and proper cycles $\sigma_i$ of length $l_i$ so that $\tau - 1 = a + 1 + \sum_{i=1}^r l_i$. It thus follows that%
\begin{align*}
  \sum_{t=0}^{\tau-2}w(\alpha(t)) &\leq n \max_{A_i \in \AC}w(A_i) + w^*_{\rm m}\sum_{i=1}^r l_i \\
	&\leq n \log_2 n + \tau w^*_{\rm m},%
\end{align*}
leading to%
\begin{equation*}
  \lim_{\tau \rightarrow \infty}\frac{1}{\tau}\max_{\alpha \in W_{\infty}(\GC)}\sum_{t=0}^{\tau-2} w(\alpha(t)) \leq w^*_{\rm m}.%
\end{equation*}
To show the converse inequality, consider an $\alpha \in W_{\infty}(\GC)$ that traces a proper cycle with mean weight equal to the maximum cycle mean $w^*_{\rm m}$. Let $l$ be the length of the cycle and write $\tau - 1 = rl + a$ for any $\tau > 1$, where $r\geq0$ and $0 \leq a < l$ are integers. This implies%
\begin{equation*}
  \frac{1}{\tau}\sum_{t=0}^{\tau-2}w(\alpha(t)) \geq \frac{lr}{\tau} w^*_{\rm m},%
\end{equation*}
and hence%
\begin{equation*}
  \frac{1}{\tau}\max_{\alpha \in W_{\infty}(\GC)}\sum_{t=0}^{\tau-2}w(\alpha(t)) \geq \frac{1}{\tau} w^*_{\rm m}(\tau - 1 - n).%
\end{equation*}
It now easily follows that $\bar{h}(\AC,G) \geq w^*_{\rm m}$, which concludes the proof. \hfill \qed
\end{pf}

\section{Relationship between the upper bounds}

In this section, we prove that in the deterministic case, where the obtained upper bound of the IE for deterministic systems and the one for uncertain ones both apply, these bounds are related by an inequality.%

Consider a deterministic system $\Sigma$ as in \eqref{eq_det_sys}, a compact controlled invariant set $Q \subseteq X$, and an invariant partition $(\AC,G)$ with $\AC = \{A_1,\ldots,A_q\}$. Let $\BC = \{B_1,\ldots,B_r\}$ be a refinement of $\AC$ and construct the weighted directed graph $\GC$ as described in Section \ref{sec_det_ub}. The sets $\WC_N(\BC,\AC)$ and $h(\BC,\AC)$ and the transition matrix $\Gamma$ are defined as in \eqref{eq_def_Gamma}. For simplicity, we assume that $\GC$ is strongly connected, in which case we know that%
\begin{equation*}
  h(\BC,\AC) = \lim_{N \rightarrow \infty}\frac{|\WC_N(\BC,\AC)|}{N} = \log_2 \rho(R),%
\end{equation*}
where $R$ is the adjacency matrix associated with a right-resolving graph.%


\begin{prop}
Given the invariant partition $(\AC,G)$, for any refinement $\BC$ of $\AC$, it holds that%
\begin{equation*}
  \bar{h}(\AC,G) \geq h(\BC,\AC) = \log_2 \rho(R).%
\end{equation*}
\end{prop}

\begin{pf}
{\colorb We use $\WC_{\tau}(\AC)$ to refer to the set $\WC_N(\BC,\AC)$, which is defined in Section~\ref{sec_det_ub}, for the case when $\BC = \AC$ and $\tau = N$, i.e., $\WC_{\tau}(\AC):= \{\alpha\in\AC^N : T_{\CC}(\alpha_i)\cap\alpha_{i+1}\neq\emptyset \ \ \forall i\in[0;N-2]\}$.}
Constructing the graph $\GC$ associated with $\AC$ as in Section \ref{sec_ub_uncertain}, leads to%
\begin{equation*}
  W_{\tau}(\GC) = \WC_{\tau}(\AC) \mbox{\quad for all\ } \tau \in \N.%
\end{equation*}
From \cite[Lem.~2]{tomar2020invariance} and \eqref{eq_thm1_ingr1}, we obtain%
\begin{equation}\label{eq_ineq}
  |\WC_{\tau}(\AC)| = |W_{\tau}(\GC)| \leq \NC(W_{\tau}(\GC)) = \bar{r}_{\inv}(\tau,Q,\AC,G).%
\end{equation}
Then \eqref{eq_ineq} yields%
\begin{equation*}
  h(\AC) := \lim_{\tau \rightarrow \infty}\frac{\log_2 |\WC_{\tau}(\AC)|}{\tau} \leq \bar{h}(\AC,G).%
\end{equation*}
It is clear that $h(\BC,\AC) \leq h(\AC)$. Hence,%
\begin{equation*}
  \log_2 \rho(R) = h(\BC,\AC) \leq h(\AC) \leq \bar{h}(\AC,G).%
\end{equation*}
This concludes the proof. \hfill \qed
\end{pf}

\section{Examples}
In this section, we illustrate the effectiveness of our proposed results on some case studies.
\subsection{A linear discrete-time system}\label{sec_lindet}

Consider the following linear control system obtained from a similarity transformation applied to the system in Example \ref{ex1}:%
\begin{equation*}
	x_{k+1} = Ax_k + \left[\begin{array}{c} 0.9463 \\ 1.051 \end{array}\right]u_k, \quad
	A = \left[\begin{array}{cc} 2 & 0.0784 \\ 0.0784 & 0.5041 \end{array}\right],%
\end{equation*}
with $x_k \in \R^2$ and $u_k \in U= [-1,1]$. Consider the set $Q$ given by the inequality%
\begin{equation*}
	\left[\begin{array}{c c}
		0.0261 & -0.4993 \\
		0.9986 & 0.0523 \\
		-0.0261 & 0.4993 \\
		-0.9986 & -0.0523
	\end{array}\right] x \leq \left[\begin{array}{c}1 \\ 1 \\ 1 \\ 1\end{array}\right],\quad x\in \R^2,%
\end{equation*}
which is compact and controlled invariant.%

To compute an upper bound on the IE of $Q$, we put $Q_X := [-1.2,1.2] \times [-2.1,2.1]$, $\eta_s := [0.04,0.08]\trn$ and $\eta_i := 0.2$.%

For the parameters $\eta_s=[0.04,0.08]\trn$, $\eta_i=0.2$, Table~\ref{tb_conju_d} lists the values of $h(\BC,\AC)$ for different selections of the coarse partition $\AC$. For the same values of $\eta_s$ and $\eta_i$, the obtained value for the bound in Theorem \ref{thm1} is $w^*_{\rm m}(\GC) = 2.5849$ (with computation time $0.048$ sec). For \texttt{dtControl} parameters Classifier $=$ `logreg' and Determinizer $=$ `maxfreq', Table~\ref{tb:conju:longerTau} shows the variation of the upper bound $h(\BC,\AC)/\tau$ with increasing control sequence length $\tau$.%

\begin{table}[]
		\centering
		\caption{Entropy estimates for Example~\ref{sec_lindet} with different choices of the determinization options in \texttt{dtControl}. Here, $h_{\inv}(Q) = 1.003$.}
		\begin{tabular}{|l|l|l|l|l|}
			\hline
			Classifier  & Determinizer & $|\AC|$ & $h(\BC,\AC)$ & time  \\ \hline\hline
			cart              & maxfreq      & 10   & 1.2133 & 10 sec \\ \hline
			logreg            & maxfreq      & 9   & 1.1802 & 10 sec \\ \hline
			linsvm            & maxfreq      & 10   & 1.2133 & 10 sec \\ \hline
			cart              & minnorm      & 135   & 1.7848 & 9 sec \\ \hline
			logreg            & minnorm      & 111   & 1.8015 & 11 sec\\ \hline 
			linsvm            & minnorm      & 143   & 1.8300 & 10 sec \\ \hline 
		\end{tabular}		
\label{tb_conju_d}
\end{table}

\begin{table}[]
	\centering
	\caption{Entropy estimates for Example~\ref{sec_lindet} with control sequences of length $\tau$. Here, $h_{\inv}(Q) = 1.003$.}
	\begin{tabular}{|l|l|l|}
		\hline
		$\tau$ & $h(\BC,\AC)/\tau$ & time     \\ \hline\hline
		1  & 1.1802 & 9.6 sec\\ \hline
		2  & 1.0688 & 16.7 sec \\  \hline
		3  & 1.0588 & 1 min 11 sec \\ \hline
	\end{tabular}
\label{tb:conju:longerTau}
\end{table}

\subsection{A scalar continuous-time nonlinear control system}\label{sec_scalarNonlinear}

Consider the following scalar continuous-time control system discussed in \cite[Ex.~7.2]{kawan2013invariance}:%
\begin{equation*}
  \Sigma:\quad \dot{x} = (-2b\sin x \cos x - \sin^2 x + \cos^2 x) + u \cos^2 x,%
\end{equation*}
where $u \in [-\rho,\rho]$, $b>0$ and $0 < \rho < b^2 + 1 =: a$. The equation describes the projectivized linearization of a controlled
damped mathematical pendulum at the unstable position, where the control acts as a reset force. The following set is controlled invariant:%
\begin{equation*}
  Q = \left[\arctan(-b-\sqrt{a + \rho}),\arctan(-b - \sqrt{a - \rho})\right].%
\end{equation*}
In fact, $Q$ is the closure of a maximal set of complete approximate controllability. With $\TC_s \in \R_{>0}$ as the sampling time, we first obtain a discrete-time system as in \eqref{eq_det_sys}. Theory suggests that the following formula holds, see\footnote{The factor $\ln(2)$ appears due to the choice of the base-$2$ logarithm
instead of the natural logarithm, which is typically used for continuous-time systems.} \cite[Ex.~7.2]{kawan2013invariance}:%
\begin{equation*}
  h_{\inv}(Q) = \frac{2}{\ln 2}\sqrt{a - \rho}.%
\end{equation*}
Discretizing the given system with sampling time $\TC_s$ results in a discrete-time system $\Sigma^{\TC_s}$ that satisfies%
\begin{equation*}
  h_{\inv}(Q;\Sigma^{\TC_s}) \geq \TC_s \cdot h_{\inv}(Q) = \frac{2\TC_s}{\ln 2} \sqrt{a - \rho}.%
\end{equation*}
The inequality is due to the fact that continuous-time open-loop control functions are lost due to the sampling (since only the piecewise constant control functions, constant on each interval of the form $[k\TC_s,(k+1)\TC_s)$, $k \in \Z_+$, are preserved under sampling). Since $Q$ can be made invariant by constant control inputs only, $Q$ is also a controlled invariant set of $\Sigma^{\TC_s}$. Table~\ref{tb:scalarNonlinear:tau:rho1b1} and \ref{tb:scalarNonlinear:tau:rho50b10} list the values of $h(\BC,\AC)/\TC_s$ for different choices of the sampling time with the parameters $(\rho = 1$, $b=1$, $\eta_s = 10^{-6}$, $\eta_i = 0.2\rho$) and ($\rho = 50$, $b = 10$, $\eta_s = 10^{-6}$, $\eta_i = 0.2\rho$), respectively. In both tables, the \verb+dtControl+ parameters are Classifier $=$ `cart' and Determinizer $=$ `maxfreq'. Table~\ref{tb:scalarNonlinear:CandD} shows the values of $h(\BC,\AC)/\TC_s$ for different selections of the coarse partition $\AC$ with the parameters $\TC_s = 0.01$, $\eta_s = 10^{-6}$, $\eta_i = 0.2\rho$, $\rho = 1$, $b=1$. For the same selection of parameters as in Table~\ref{tb:scalarNonlinear:tau:rho1b1} with $\TC_s = 0.01$, Table~\ref{tb:scalarNonlinear:ulength} presents the variation of the upper bound $h(\BC,\AC)/(\tau\TC_s)$ with increasing length $\tau$ of the control sequences.%

\begin{table}[]
	\centering
	\caption{Entropy estimates for Example~\ref{sec_scalarNonlinear} with $\rho=1$, $b=1$ and different choices of the sampling time $\TC_s$. Here, $h_{\inv}(Q) = 2.8854$.}
	\begin{tabular}{|l|l|l|l|}
		\hline
		$\TC_s$ & $|\AC|$ & $h(\BC,\AC)/\TC_s$ & time  \\  \hline\hline
		0.8    & 11            & 4.0207  & 21.23 hr  \\ \hline
		0.5    & 6             & 4.0847  & 2.98 hr    \\ \hline
		0.1    & 2             & 4.744   & 3.33 min     \\ \hline
		0.01   & 2             & 5.1994  & 55 sec    \\ \hline
		0.001  & 2             & 24.7    & 60 sec     \\ \hline
	\end{tabular}
	\label{tb:scalarNonlinear:tau:rho1b1}
\end{table}

\begin{table}[]
	\centering
	\caption{Entropy estimates for Example~\ref{sec_scalarNonlinear} with $\rho=50$, $b=10$ and different choices of the sampling time $\TC_s$. Here, $h_{\inv}(Q) = 20.6058$.}
	\begin{tabular}{|l|l|l|l|}
		\hline
		$\TC_s$ & $|\AC|$ & $h(\BC,\AC)/\TC_s$ & time \\ \hline\hline
		0.11   & 15            & 28.5012  & 1.9 hr    \\ \hline
		0.1    & 11            & 29.1723  &  1.35 hr    \\ \hline
		0.01   & 2             & 34.4707 & 13 sec     \\ \hline
		0.001  & 2             & 55.5067 & 12 sec       \\ \hline
		0.0001 & 2             & 1.5635e+03 & 31 sec    \\ \hline
		
	\end{tabular}
	\label{tb:scalarNonlinear:tau:rho50b10}
\end{table}

\begin{table}[]
	\centering
	\caption{Entropy estimates for Example~\ref{sec_scalarNonlinear} with different choices of \texttt{dtControl} parameters. Here, $h_{\inv}(Q) = 2.8854$.}
	\begin{tabular}{|l|l|l|l|l|}
		\hline
		Classifier & Determinizer & $|\AC|$ &   $h(\BC,\AC)/\TC_s$   & time  \\ \hline\hline
		cart       & maxfreq      & 2   & 5.1994  & 55 sec   \\ \hline
		logreg     & maxfreq      & 2   & 5.1994  & 65 sec   \\ \hline
		linsvm     & maxfreq      & 2   & 5.1994  & 61 sec   \\ \hline
		cart       & minnorm      & 11  & 6.4475  & 57 sec   \\ \hline
		logreg     & minnorm      & 11  & 6.4475  & 74 sec   \\ \hline
	\end{tabular}
	\label{tb:scalarNonlinear:CandD}
\end{table}

\begin{table}
	\centering
	\caption{Upper bound $h(\BC,\AC)/(\tau\TC_s)$ for Example~\ref{sec_scalarNonlinear} with control sequences of length $\tau$, Classifier = `cart', and Determinizer = `maxfreq' in \texttt{dtControl}. Here, $h_{\inv}(Q) = 2.8854$.}
	\begin{tabular}{|l|l|l|}
		\hline
		$\tau$ & $h(\BC,\AC)/(\tau\TC_s)$  & time   \\ \hline\hline
		1  & 5.1994 & 57 sec \\ \hline
		2  & 5.0036 & 7.5 min \\ \hline
		3  & 4.9547 & 1.91 hr \\ \hline
		4  & 4.9266 & 27.27 hr \\ \hline
	\end{tabular}
	\label{tb:scalarNonlinear:ulength}
\end{table}

\subsection{A 2d uniformly hyperbolic set}\label{sec_hyperbolic}

Consider the map%
\begin{equation*}
  f(x,y) := (5 - 0.3y - x^2,x),\quad f:\R^2 \rightarrow \R^2,%
\end{equation*}
from the H\'enon family, one of the most-studied classes of dynamical systems that exhibit chaotic behavior. We extend $f$ to a control system with additive control:%
\begin{equation*}
	\Sigma:\quad \left[\begin{array}{c} x_{t+1} \\ y_{t+1} \end{array}\right] = \left[\begin{array}{c} 5 - 0.3 y_t - x_t^2 + u_t \\ x_t + v_t \end{array}\right],%
\end{equation*}
with $\max\{|u_t|,|v_t|\} \leq \ep$. It is known that $f$ has a non-attracting uniformly hyperbolic set $\Lambda$, which is a topological horseshoe. This set is contained in the square centered at the origin with side length \cite[Thm.~4.2]{robinson1998dynamical}%
\begin{equation*}
  r := 1.3 + \sqrt{(1.3)^2 + 20} \approx 5.9573.%
\end{equation*}
If the size $\ep$ of the control range is chosen small enough, the set $\Lambda$ is ``blown up'' to a compact controlled invariant set $Q^{\ep}$ with nonempty interior which is not much larger than $\Lambda$, see \cite[Sec.~6]{kawan2020control}. Moreover, the theory suggests that as $\ep\downarrow0$, $h_{\inv}(Q^{\ep})$ converges to the negative topological pressure of $f|_{\Lambda}$ w.r.t.~the negative unstable log-determinant on $\Lambda$; see \cite{bowen2008equilibrium} for definitions. A numerical estimate for this quantity, obtained in \cite[Table 2]{froyland1999using} via Ulam's method, is $0.696$.%

We select $\tilde{Q} = [-r/2,r/2]^2$. For $\ep = 0.08$, using \verb+SCOTS+ with parameter values $\eta_s = [0.009,0.009]\trn$ and $\eta_i = [0.01,0.01]\trn$, through iteration, we obtain an all-time controlled invariant set $Q \subset \tilde{Q}$. In the iteration, we begin with the set $\tilde{Q}$ and, as the first step, we compute an invariant controller for the system $\Sigma$. Let $Q_1$ be the domain of the obtained controller. Consider the time-reversed system%
\begin{equation*}
	\Sigma^-:\quad \left[\begin{array}{c} x_{t+1} \\ y_{t+1} \end{array}\right] = \left[\begin{array}{c}  y_t - v_t \\ \frac{1}{0.3}(5 - y_t^2 + u_t - x_t) \end{array}\right].%
\end{equation*}
In the second step, we compute an invariant controller for $\Sigma^-$ in the set $Q_1$, and denote the controller domain by $Q_2$. In the third step, we compute an invariant controller for $\Sigma$, but in the set $Q_2$, and denote the controller domain by $Q_3$. The steps are repeated until $Q_i = Q_{i+1} =: Q$. In this way, we hope to approximate $Q^{\ep}$.%

\begin{figure}
\centering
\includegraphics[width=0.4\textwidth]{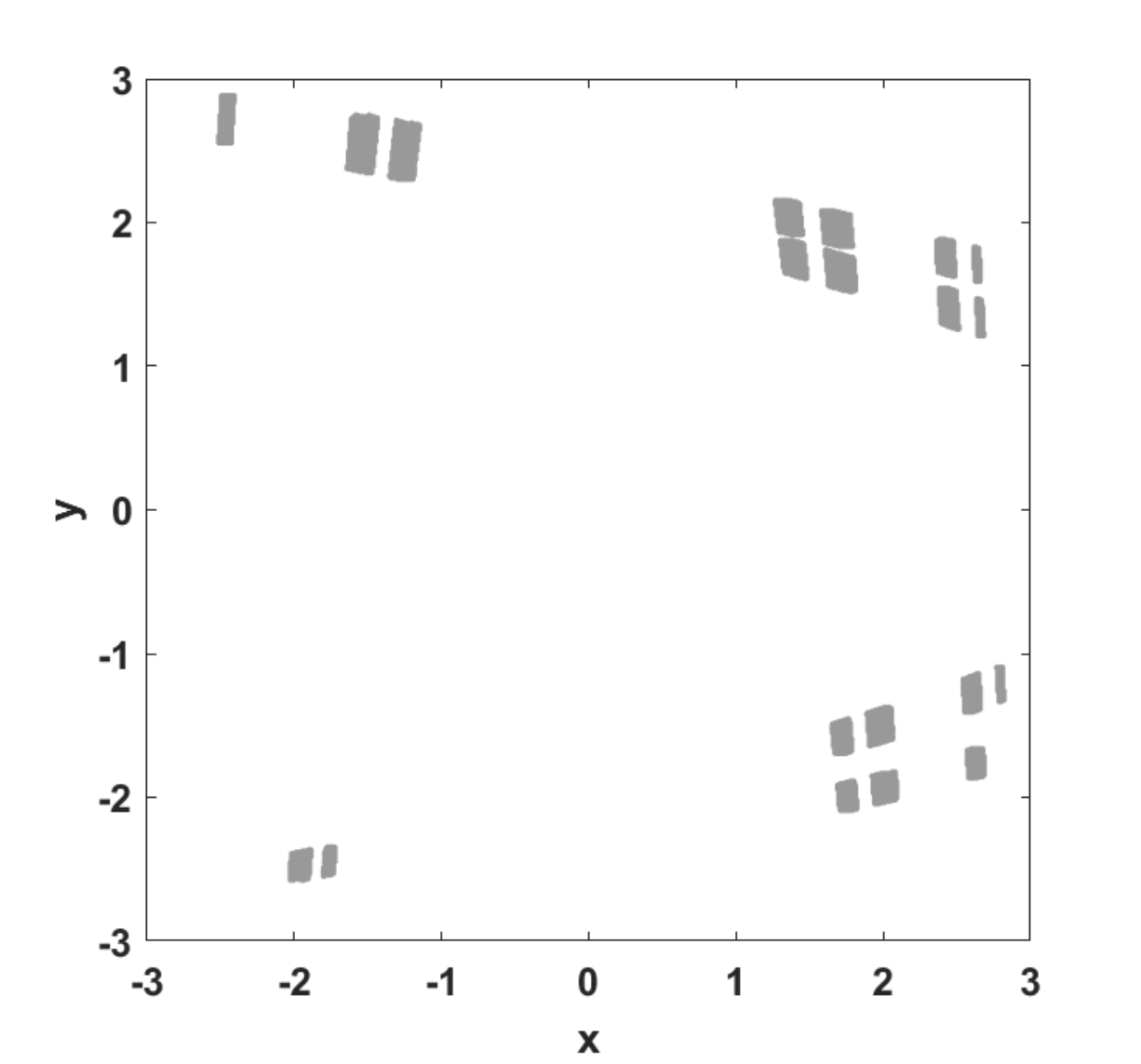}
\caption{The set $Q$ for Example~\ref{sec_hyperbolic}}.
\label{fg:hyperbolic:domain}
\end{figure}

Figure \ref{fg:hyperbolic:domain} shows the set $Q$. For the parameter values $\varepsilon=0.08$, $\eta_s = [0.009,0.009]\trn$, $\eta_i = [0.01,0.01]\trn$, Table \ref{tb:hyperbolic} lists the values of $h(\BC,\AC)$ for different choices of the coarse partition $\AC$. For the same values of $\varepsilon$, $\eta_s$ and $\eta_i$, the obtained value for the bound in Theorem \ref{thm1} is $w^*_{\rm m}(\GC) = 3.5646$ (with computation time $2.51$ sec).%

\begin{table}[]
\centering
\caption{Entropy estimates for Example~\ref{sec_hyperbolic} with different selections of \texttt{dtControl} options. Here, $h_{\inv}(Q) \approx 0.696$.}
\begin{tabular}{|l|l|l|l|l|}
	\hline
	Classifier & Determinizer & $|\AC|$   &   $h(\BC,\AC)$ & time     \\ \hline\hline
	cart       & maxfreq      & 573  & 2.3884 & 0.95 min \\ \hline
	linsvm     & maxfreq      & 567  & 2.3956 & 1.82 min \\ \hline
	logreg	   & maxfreq	  & 454  & 2.3994 & 1.4 min \\ \hline
	cart       & minnorm      & 1921  & 2.9342 & 1 min\\ \hline
	logreg     & minnorm      & 1533  & 2.9215 & 2 min \\ \hline
	linsvm     & minnorm      & 1923  & 2.9376 & 2.15 min \\ \hline
\end{tabular}
\label{tb:hyperbolic}
\end{table}

\subsection{An uncertain linear system}\label{sec_ulinear}

\begin{figure}
\centering
\includegraphics[scale=.4]{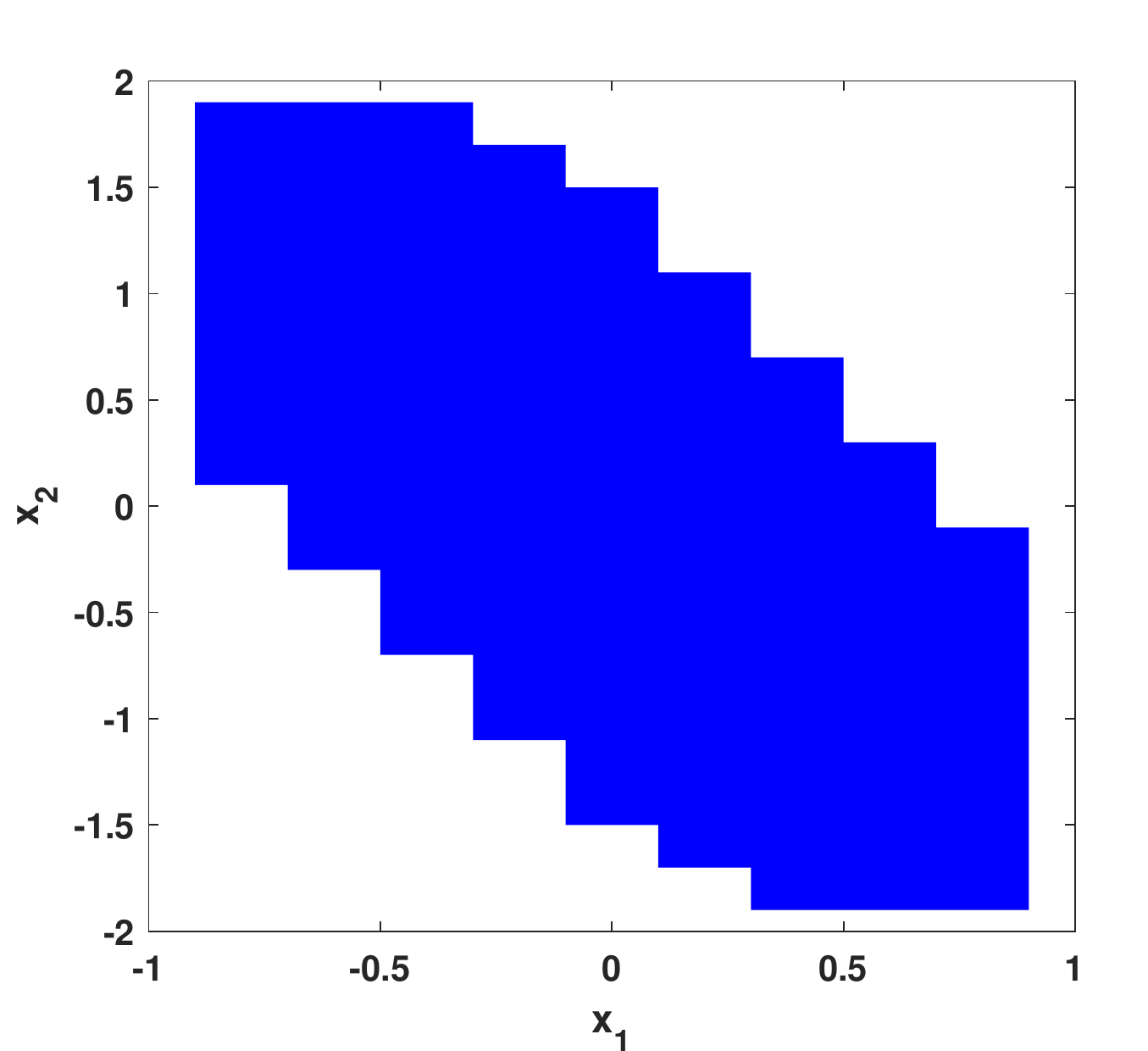}
\caption{The set $Q$ in Example~\ref{sec_ulinear} which is defined as the domain of the invariant controller computed from \texttt{SCOTS}.}
\label{fg:ex:ulinear:Q}
\end{figure}

We consider an uncertain linear control system%
\begin{align*}
  \Sigma:\quad x_{t+1} \in \left[\begin{array}{cc} 2 & 1 \\ -0.4 & 0.5 \end{array}\right]x_t + \left[\begin{array}{c} 1 \\ 1 \end{array}\right]u_t + W%
\end{align*}
with $x_t \in \R^2$, $u_t \in U := [-1,1]$, and the disturbance set $W := [-0.1,0.1]^2$. For a set $Q \subseteq [-1,1] \tm [-2,2]$, we compute an upper and a lower bound of $\bar{h}_{\inv}(Q)$. We used \verb+SCOTS+ to obtain an invariant controller for the state set $[-1,1] \tm [-2,2]$ with $[0.2,0.2]\trn$ and $0.05$ as the state and input set grid parameters, respectively. The set $Q$ is taken to be the domain of the obtained controller that consists of 109 state grid cells each of size $0.2 \tm 0.2$. Figure~\ref{fg:ex:ulinear:Q} shows the set $Q$.%

{\bf Computation of the lower bound:} We utilize \cite[Thm.~7]{tomar2020invariance} to compute a lower bound. From \cite[Rem.~2]{tomar2020invariance}, we know that the lower bound in \cite[Thm.~7]{tomar2020invariance} is invariant under coordinate transformations. After a similarity transformation which diagonalizes the dynamical matrix, we have%
\begin{equation*}
  z_{t+1} = \left[\begin{array}{cc} 1.6531 & 0 \\ 0 & 0.8469 \end{array}\right]z_t + V^{-1}\left[\begin{array}{c} 1 \\ 1 \end{array}\right] + V^{-1}W,%
\end{equation*}
where $V=\left[\begin{array}{cc}0.9448 & -0.6552\\ -0.3277 & 0.7555\end{array}\right]$.
For $i = 1,2$, let $\pi_i$ denote the canonical projection to the $i$th coordinate. Then, $\pi_1 V^{-1}Q = [-2.1207,2.1207]$, $\pi_2V^{-1}Q = [-3.4,3.4]$, $\pi_1V^{-1}W = [-0.2827,0.2827]$ and $\pi_2V^{-1}W = [-0.2550,0.2550]$. By \cite[Thm.~7]{tomar2020invariance}, one obtains%
\begin{equation*}
   0.9316 \leq \bar{h}_{\inv}(Q).%
\end{equation*}

{\bf Computation of the upper bound:} We construct an invariant partition $(\AC,G)$ of $Q$ by selecting the set of grid cells in the domain of the controller obtained from \verb+SCOTS+ as the cover $\AC$. Let $C:\AC \rightrightarrows U$ denote the controller from \verb+SCOTS+. For $A \in \AC$, $C(A)$ is the list of control inputs in the controller assigned to cell $A$ such that each of the control inputs in the list ensures invariance of the states in $A$ w.r.t.~the set $Q$. For each $A \in \AC$, we define $G(A) := u \in C(A)$, where $u$ is chosen such that $F(A,u)$ has nonempty intersection with a minimum number of elements of $\AC$. If there are multiple such control values, then one of them is selected randomly. Using $(\AC,G)$ and the transition function $F$ of the system, we construct a weighted directed graph $\GC$ as described in Section \ref{sec_ub_uncertain}. We used the \verb+LEMON library+ to compute the maximum cycle mean weight for the graph $\GC$ and obtained $w^*_{\rm m}(\GC) = 3.3219$ with computation time 0.027 sec. Thus, $\bar{h}_{\inv}(Q) \leq 3.3219$.%

\textbf{Discussion on the selection of partition:} A better upper bound is expected when the number of outgoing edges, for every node in the graph, is smaller. As a heuristic, gradually smaller values of the state grid parameter $\eta_s$ can be tried. But very small $\eta_s$ that make width of the grid cell smaller than that of the disturbance set should be avoided, because in that case, the number of outgoing edges for any cell will begin to rise. This can also be observed from Table \ref{tbl_uncertain}.%

\begin{table}[]
	\centering
	\caption{Entropy estimates for Example~\ref{sec_ulinear}, with $\eta_i = 0.05$.}
	\begin{tabular}{|c|c|c|}
		\hline
		$\eta_s$ & $w_{\rm m}^*$  & time(sec) \\
		\hline \hline
		0.03 & $6.4594 $ & $1.112$\\ \hline
		0.06 & $5 $ & $0.129$\\ \hline
		0.09 & $4.2811 $ & $0.051$\\ \hline
		0.1 & $4.3923 $ & $0.033$\\ \hline
		0.2 & $3.3219 $ & $0.027$\\ \hline
	\end{tabular}
	\label{tbl_uncertain}
\end{table}

\section{Software tools and pseudo-code}

In this section, we provide brief descriptions of the used software tools and summarize our algorithms in terms of pseudo-code (cf. Figures \ref{proced:d} and \ref{proced:u}).%

\textbf{Description of the computation of the maximum cycle mean (MCM) using the LEMON library:} The maximum cycle mean of a directed weighted graph can be computed by Karp's algorithm which runs in $O(nm)$ time, where $n$ and $m$ are the number of nodes and edges in the graph, respectively. For MCM, we utilize LEMON\footnote{\url{https://lemon.cs.elte.hu/trac/lemon}} which is a \texttt{C++} library that provides efficient implementations of algorithms related to graphs. LEMON provides the implementation of Karp's algorithm in the class \texttt{KarpMmc}. The class constructor requires two arguments: \texttt{Digraph} and \texttt{CostMap}. \texttt{Digraph} specifies the type of directed-graph implementation to be used, while \texttt{CostMap} is a map that specifies the weight assigned to each edge in the graph. The algorithm can be executed by the class member function \texttt{findCycleMean()}. Although the implementation computes the minimum cycle mean, the maximum cycle mean can be easily computed by assigning a negative sign to every edge weight.%

\textbf{Description of dtControl}: \texttt{dtControl} is a software tool, written in Python, for post-processing memoryless symbolic controllers into various compact and more interpretable representations~\cite{ashok2020dtcontrol}. It provides parameters like `Classifier' and `Determinizer' to adjust the DT (decision tree) learning algorithm. The classifier \texttt{cart} (classification and regression tree) allows only axis-aligned splits of the state space at any node, while the linear classifiers as \texttt{logreg} (logistic regression) and \texttt{linsvm} (linear support vector machine) allow oblique splits. With axis-aligned splits, the elements of the coarse partition are hyperrectangles, while with oblique splits more general partition elements, in the form of union of hyperrectangles can be obtained. With determinizer \texttt{minnorm}, the controller is first determinized by the selection of control values with the smallest norm and then the decision tree is learned. The determinizer \texttt{maxfreq} gives the best upper bounds in all examples. Let $C:X\rightrightarrows U$ denote the controller that is fed to \texttt{dtControl} and $S_n$ denote the subset of the state space corresponding to the node $n$ of the decision tree. When \texttt{maxfreq} is selected, then, during the construction of the decision tree, for every node $n$, the corresponding part of the controller $(C|_{S_n})$ is determinized through the selection of the control values that have the maximum frequency of appearance in the set $\cup_{x\in S_n}C(x)$. Then a classifier is learned for the subset $S_n$. \texttt{Maxfreq} typically leads to quite small decision trees, thus coarser partitions.%

\begin{algorithm}
\captionof{figure}{Procedure for computing an upper bound of IE in the deterministic case}
\label{proced:d}
\hspace*{\algorithmicindent} \textbf{Input:} $\eta_s$, $\eta_i$, $Q_X$, $Q$, $U$, $f$, $\tau$, classifier, determinizer \\
    \hspace*{\algorithmicindent} \textbf{Output:} $h(\BC,\AC)$
\begin{algorithmic}[1]
\State $C$ $\gets$ \texttt{SCOTS}($\eta_s$, $\eta_i$, $Q_X$, $Q$, $U$, $f$, $\tau$)
\State $(\bar{C},\AC) \gets $ \texttt{dtControl}(classifier, determinizer, $C$)
\State $\GC \gets$ graph($f$, $\bar{C}$, $\AC$)
\State $\GC_{\textrm{scc}} \gets $ stronglyConnectedComponents($\GC$)
\State $l \gets \{\}$
\ForAll{$\GC_k \in \GC_{\textrm{scc}}$} 
\State  $\bar{\GC}_k \gets$ rightResolvingGraph($\GC_k$)
\State $R^k \gets $ adjacencyMatrix($\bar{\GC}_k $)
\State $l \gets l \cup \{\log_2\rho(R^k) \}$
\EndFor
\State $h(\BC,\AC) \gets \max(l)$
\end{algorithmic}
\end{algorithm}

In Fig.~\ref{proced:d}, $C$ is the set valued map $C:\BC\rightrightarrows U^{\tau}$obtained from \texttt{SCOTS} while $\bar{C}$ is a single valued map $\bar{C}:\BC\to U^{\tau}$. \texttt{SCOTS} runs in $O(n^2m)$ time, where $n$ and $m$ are the number of cells in the state and the input grid, respectively. The decision tree learning takes $O(n^2\log(n))$ time.\footnote{\url{https://scikit-learn.org/stable/modules/tree.html\#complexity}} Each element of the coarse partition $\AC$ corresponds to some leaf node in the decision tree obtained from \texttt{dtControl}. $\GC_{\textrm{scc}}$ is the set of strongly connected components of the graph $\GC$, and it can be computed in $O(n+E)$ time, where $n$ and $E$ are the number of nodes and edges, respectively, in the graph. A right-resolving graph can be computed by the power-set construction in $O(2^n|\AC|)$. The adjacency matrix and the right-resolving graph are computed simultaneously. {Thus, the procedure in Fig.~\ref{proced:d} runs in $O(n2^n)$.}%

\begin{algorithm}
\captionof{figure}{Procedure for computing an upper bound of IE in the uncertain case}
\label{proced:u}
\hspace*{\algorithmicindent} \textbf{Input:} $\eta_s$, $\eta_i$, $Q_X$, $Q$, $U$, $F$ \\
    \hspace*{\algorithmicindent} \textbf{Output:} $w_m^*(\GC)$
\begin{algorithmic}[1]
\State $C$ $\gets$ \texttt{SCOTS}($\eta_s$, $\eta_i$, $Q_X$, $Q$, $U$, $F$, $1$)
\State $\bar{C} \gets $ determinize($C$)
\State $\GC \gets$ graph($F$, $\bar{C}$)
\State $w_{\rm m}^*(\GC) \gets $ maxCycleMean($\GC$)
\end{algorithmic}
\end{algorithm}

In Fig.~\ref{proced:u}, the controller $C:\BC\rightrightarrows U$ can be determinized through the selection of such control inputs that result in minimum number of successor state-cells, in $O(\bar{n}m)$ time, where $\bar{n}=|\BC|$. The maximum cycle mean by Karp's algorithm can be computed in $O(\bar{n}\bar{m})$, where $\bar{m}$ is the number of edges in the graph, which in the worst case will be $\bar{n}^2$. Thus, the procedure in Fig.~\ref{proced:u} runs in $O(n^2m+n^3)$.%



{\colorb
\begin{remark}
    To reduce computational complexity, for the uncertain case, one can leverage the proposed compositionality results in~\cite{tomar2020compositional} for the computation of an overapproximation of the invariance entropy for a large scale interconnected system in a divide and conquer manner by computing overapproximations for subsystems using the method proposed here. Thus complexity breaks down to the level of subsystems.
\end{remark}
}

The code is publicly accessible at \url{https://github.com/
mahendrasinghtomar/Invariance_Entropy_
upper_bounds}.

{\bf Quality of upper bounds of IE:}
For uncertain nonlinear systems, because of the absence of any theory providing a lower bound (better than zero) in the literature, we do not know how far our computed upper bounds are from the actual values. For uncertain linear systems with additive disturbance, one can comment on this gap based on the availability of a lower bound~\cite[Thm. 7]{tomar2020invariance}. In the deterministic nonlinear case, the gap is not yet quantified as well.%

\section{Conclusion and future work}

Our first contribution is the combination of three different algorithms designed for different purposes to numerically compute an upper bound of the invariance entropy of deterministic control systems. The second contribution is a procedure to numerically compute an upper bound for the invariance entropy of uncertain control systems. We also describe the relationship between the two upper bounds and thus the need for the second bound. Finally, we illustrate the effectiveness of the proposed procedures on four examples. Open questions for future work include the selection of entropy-minimizing partitions and the computation of lower bounds of IE for uncertain nonlinear systems.%



\end{document}